\documentclass[a4paper,12pt]{article}
\usepackage{amsmath,amssymb,amsfonts,amsthm,eucal}
\usepackage{graphicx}
\usepackage{tikz}
\usetikzlibrary{positioning,decorations.pathreplacing}
\usepackage{float}
\usepackage{booktabs,float,rotating}

\newtheorem{proposition}{Proposition}

\newtheorem{theorem}{Theorem}
\newtheorem{definition}{Definition}

\newtheorem{example}{Example}

\newtheorem{axiom}{Axiom}

\begin{document}
\date{\today}
\title{On priority in multi-issue bankruptcy problems with crossed claims}

\author{Rick K. Acosta-Vega\thanks{Faculty of Engineering, University of Magdalena, Colombia. \{racosta@unimagdalena.edu.co\}}
\and Encarnaci\'{o}n Algaba\thanks{%
Department of Applied Mathematics II and IMUS, University of Seville, Spain. \{ealgaba@us.es\}}
\and Joaqu\'in S\'anchez-Soriano\thanks{%
\textbf{Corresponding author}. R.I. Center of Operations Research (CIO), Miguel Hern\'andez University of Elche, Spain. \{joaquin@umh.es\}}
}
\maketitle

\begin{abstract}
In this paper, we analyze the problem of how to adapt the concept of priority to situations where several perfectly divisible resources have to be allocated among certain set of agents that have exactly one claim which is used for all resources. In particular, we introduce constrained sequential priority rules and two constrained random arrival rules, which extend the classical sequential priority rules and the random arrival rule to these situations. Moreover, we provide an axiomatic analysis of these rules. 
\end{abstract}

\section{Introduction}\label{intro}

There is a vast literature on problems in which a resource must be allocated among a set of claimants. A relevant allocation problem arises when there is a (perfectly divisible) resource (for example, money, water, time...) over which there is a set of claimants who have rights or demands, but the resource is scarce. This problem is known as bankruptcy problem (O'Neill, 1982; Aumann and Maschler, 1985). Many allocation rules have been defined to provide a solution to this problem (see Thomson (2003, 2015, 2019) for a detailed inventory of rules and their axiomatic analysis). One of the most prominent allocation rules is the random arrival rule (O'Neill, 1982) also known as the run-to-the-bank rule (Young, 1994). For this rule the following dynamic interpretation can be given. The problem is solved by holding a race in such a way that the first claimant who arrives is satisfied with as much of the resource, the second is satisfied with as much of the resource as is left, and so on until the resource is exhausted. If we consider that all orders of arrival are equally likely, the random arrival rule is the expected value. Likewise, we can consider that each of these orders of arrival is a priority relationship defined ex ante over the claimants according to some criterion, and thus we would have a rule for each of these orders, these rules are called sequential priority rules. Moreover, S\'anchez-Soriano (2021) extends the sequential priority rules and, their average, the random arrival rule to the case in which withdrawals are bounded from above.

An extension of bankruptcy problems are multi-issue allocation problems (Calleja et al., 2005). In these situations there is also a (perfectly divisible) resource but it can be distributed between several issues, and a (finite) number of agents have claims on each of those issues, so that the total claim is higher than the available resource. This problem is also solved by means of allocation rules and there are different ways to do it. Calleja et al. (2005) introduce two extensions of the random arrival rule to the context of multi-issue allocation problems. Gonz\'alez-Alc\'on et al. (2007) propose a new extension of the random arrival rule to this context. An alternative approach to multi-issue bankruptcy problems is when the resource has already been divided a priori into the issues according to exogenous criteria but claimants continue having a claim for each issue (Izquierdo and Timoner, 2016).

\begin{figure}
\begin{center}
\begin{tikzpicture}
\node[draw] (c8) at (10.5,0) {{\small C8}};
\node[draw] (c7) at (9,0) {{\small C7}};
\node[draw] (c6) at (7.5,0) {{\small C6}};
\node[draw] (c9) at (6,0) {{\small C5}};
\node[draw] (c5) at (4.5,0) {{\small C4}};
\node[draw] (c4) at (3,0) {{\small C3}};
\node[draw] (c3) at (1.5,0) {{\small C2}};
\node[draw] (c2) at (0,0) {{\small C1}};

\node[draw] (e1) at (1.5,2.5) {{\small E1}};
\node[draw] (e2) at (5,2.5) {{\small E2}};
\node[draw] (e3) at (8,2.5) {{\small E3}};

\draw [->] (c2) -- (e1);
\draw [->] (c3) -- (e1);
\draw [->] (c4) -- (e1);
\draw [->] (c5) -- (e1);

\draw [->] (c4) -- (e2);
\draw [->] (c5) -- (e2);
\draw [->] (c6) -- (e2);
\draw [->] (c7) -- (e2);
\draw [->] (c9) -- (e2);

\draw [->] (c6) -- (e3);
\draw [->] (c7) -- (e3);
\draw [->] (c8) -- (e3);

\end{tikzpicture}
\end{center}
\caption{Multi-issue problem with crossed claims.}\label{2levels}
\end{figure}
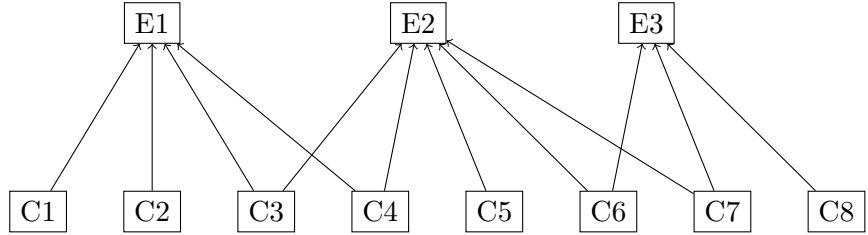

However, the situation described in Figure \ref{2levels} does not fit to any multi-issue allocation bankruptcy problem as referred in the previous paragraph, but to a multi-issue bankruptcy problem with crossed claims introduced by Acosta-Vega et al. (2021, 2022). These describe situations in which there are several (perfect divisible) resources and a (finite) set of agents who have claims on them, but only one claim (not a claim for each resource) with which one or more resources are requested. The total claim for each resource exceeds its availability.

In this paper, in order to solve allocation problems as the described in Figure \ref{2levels}, we introduce sequential priority rules and, their average, the random arrival rule for multi-issue allocation problems with crossed claims that naturally extend the sequential priority rule and the random arrival rule for single issue allocation problems.

The rest of the paper is organized as follows. Section \ref{model} presents multi-issue allocation problems with crossed claims (MBC), and the concept of rule for these problems. In Section \ref{prop}, constrained sequential priority rules and, their average, the constrained random arrival rule for multi-issue allocation problems with crossed claims are defined. In Section \ref{properties}, we present several properties which are interesting in the context of MBC problems. Section \ref{conc} concludes.


\section{Multi-issue bankruptcy problems with crossed claims}
\label{model}
Bankruptcy problems concern allocating the available resource (estate) to a number of individuals (claimants) according to their demands/claims when all these demands add up to more than the available resource. Mathematically, a bankruptcy problem is given by a triplet $(N,E,c)$, where $N$ is the set of claimants, $E \in \mathbb{R}_+$ is a perfectly divisible amount of resource (the \emph{issue} or \emph{estate}) to be divided, and $c$ is the vector of demands, such that $C = \sum_{j \in N} c_j> E$. The set of all bankruptcy problems with set of agents $N$ is denoted by $BP^{N}$.

The main goal in a bankruptcy problem is to find an allocation which is as fair as possible, taking into account the demands of the claimants. An allocation rule $R$ is a function defined for each $N$ from $BP^{N}$
to $\mathbb{R}_{+}^{N}$ such that for each $(N,E,c)\in BP^{N}$,  $R(N,E,c)\in\mathbb{R}_{+}^{N}$ satisfies the following two conditions:
\begin{enumerate}
\item $\sum_{i\in N}R_{i}(N,E,c)=E$.
\item $0\leq R_{i}(N,E,c)\leq c_{i},\:\forall i\in N$.
\end{enumerate}
Many allocation rules have been proposed depending on the principle(s) of fairness or rationale used. Allocation rules defined on the basis of the existence of a priority order on the claimants and that satisfies their demands sequentially according to that priority order are the following:

The \emph{sequential priority} rule associated with $\sigma \in \Sigma(N)$ ($SP^{\sigma}$) is defined by setting
$$
SP_{i}^{\sigma}(E,c)=min\left\{ c_{i},\:max\{0,\:E-\sum_{j\in N:\sigma(j)<\sigma(i)}c_{j}\}\right\} ,\:\forall i\in N,
$$
where $\Sigma(N)$ is the set of all possible orders of $N$. The family of sequential priority rules, SP-family, consists of all these rules.

The \emph{random arrival} rule (RA) (O'Neill, 1982) is defined by setting
$$
RA_{i}(E,c)=\frac{1}{n!}\sum_{\sigma\in\Sigma(N)}SP_{i}^{\sigma}(E,c) ,\:\forall i\in N.
$$

The random arrival rule coincides with the Shapley value (Shapley, 1953) of the pessimistic bankruptcy game (O'Neill, 1982). This is a relevant fact, because the Shapley value is one of the most prominent solution concepts in cooperative game theory (Roth, 1988; Algaba et al., 2019a, 2019b).

We now consider a situation where there is a finite set of issues  $\mathcal{I}=\{1,2,...,m\}$ such that  there is a perfectly divisible amount $e_i$ of each issue $i$. Let $E =(e_1,e_2,...,e_m)$ be the vector of available amounts of issues. There is a finite set of claimants $N=\{1,2,...,n\}$ such that each claimant $j$ claims $c_j$. Let $c=(c_{1},c_{2},...,c_{n})$ be the vector of claims. However, each claimant claims to different sets of issues, in general. Thus, let $\alpha$ be a set-valued function that associates with every $j \in N$ a set $\alpha(j) \subset \mathcal{I}$. In fact, $\alpha(j)$ represents the issues to which claimant $j$ asks for. Furthermore, $\sum_{j:i \in \alpha(j)} c_j > e_i$, for all $i \in \mathcal{I}$, otherwise, those issues could be discarded from the problem because they do not impose any limitation, and so the allocation would be trivial. Therefore, a multi-issue bankruptcy problem with crossed claims (MBC in short) is defined by a 5-tuple $(\mathcal{I},N,E,c,\alpha)$, and the family of all these problems is denoted by $\mathcal{MBC}$.

Given a problem $(\mathcal{I},N,E,c,\alpha) \in \mathcal{MBC}$, a \emph{feasible allocation} for it, it is a vector $x \in \mathbb{R}^N$ such that:
\begin{enumerate}
\item $0 \leq x_i \leq c_i$, for all $i \in N$.
\item $\sum_{i \in N:j \in \alpha(i)}x_i \leq e_j$, for all $j \in M$,
\end{enumerate}
and we denote by $A(\mathcal{I},N,E,c,\alpha)$ the set of all its feasible allocations.

An \emph{ allocation rule} for multi-issue bankruptcy problems with crossed claims is a mapping $R$ that associates with every $(\mathcal{I},N,E,c,\alpha) \in \mathcal{MBC}$ a unique feasible allocation $R(\mathcal{I},N,E,c,\alpha) \in A(\mathcal{I},N,E,c,\alpha)$.


\section{Constrained sequential priority rules for $\mathcal{MBC}$ problems}\label{prop}

Sequential priority rules are defined when there is a priority order defined over the set of claimants, in such a way that if a claimant has a higher priority than another, the first must be satisfied first in her demand. This rule simply allocates the resource according to the scheme of first come first served, where the order of arrival is given by the priority relationship. The question in $\mathcal{MBC}$ problems is what ``be satisfied with as much of the resource as is left'' means. In the context of one-issue allocation problems, ``be satisfied with as much of the resource as is left'' is a simple idea since there is only one issue. How to extrapolate this to the MBC situations. To answer this question, we introduce the constrained sequential priority rule (CSP in short) which follows the same process as sequential priority rules but taking into account that there are several resources or issues. This rule is formally defined below.

\begin{definition}\label{CSP}
Let $(\mathcal{I},N,E,c,\alpha) \in \mathcal{MBC}$, and $\sigma \in \Sigma(N)$, the \emph{constrained sequential priority rule} associated with $\sigma \in \Sigma(N)$ for $(\mathcal{I},N,E,c,\alpha)$, $CSP^\sigma(\mathcal{I},N,E,c,\alpha)$,  is defined as follows:

$$
CSP_{j}^{\sigma}(\mathcal{I},N,E,c,\alpha)=min\left\{ c_{j},\:max\left\{0,\:\min_{i \in \alpha(j)}\{e_i-\sum_{\substack{k \in N:i \in \alpha(k) \\[0.1cm]\sigma(k)<\sigma(j)}}c_{k}\}\right\}\right\} ,\:\forall j\in N,
$$
where $\Sigma(N)$ is the set of all possible orders of $N$.
\end{definition}

For each $\sigma \in \Sigma(N)$, the iterative procedure of $CSP^{\sigma}$ is well-defined and always leads to a single point. Moreover, by definition, it ends in a finite number of steps, at most $|N|$. Finally, when we have a one-issue allocation problem, then we obtain $SP^{\sigma}$. Therefore, this definition extends sequential priority rules to the context of MBC. The following example illustrates how CSP works.

\begin{example}\label{ex1}
Consider the situation described by Figure \ref{2levels} with $\mathcal{I}=\{1,2,3\}$, $N=\{1,2,3,4,5,6,7,8\}$, $E=(9,12,9)$, $c=(3,5,4,3,5,4,3,5)$, and $\alpha(1)=\{1\},\alpha(2)=\{1\},\alpha(3)=\{1,2\},\alpha(4)=\{1,2\},\alpha(5)=\{2\},\alpha(6)=\{2,3\},\alpha(7)=\{2,3\},\alpha(8)=\{3\}$. If we take for example the priority order $\sigma = 1\:3\:5\:7\:2\:4\:6\:8$, $CSP^{\sigma}$ is calculated sequentially as follows:

\bigskip

\noindent First, claimant 1 is attended:
$$
CSP_{1}^{\sigma}(\mathcal{I},N,E,c,\alpha) = 3.
$$
Second, the resources are updated down $E=(6,12,9)$ and claimant 3 is attended:
$$
CSP_{3}^{\sigma}(\mathcal{I},N,E,c,\alpha) = 4.
$$
Third, the resources are updated down $E=(2,8,9)$ and claimant 5 is attended:
$$
CSP_{5}^{\sigma}(\mathcal{I},N,E,c,\alpha) = 5.
$$
Fourth, the resources are updated down $E=(2,3,9)$ and claimant 7 is attended:
$$
CSP_{7}^{\sigma}(\mathcal{I},N,E,c,\alpha) = 3.
$$
Fifth, the resources are updated down $E=(2,0,6)$ and claimant 2 is attended:
$$
CSP_{2}^{\sigma}(\mathcal{I},N,E,c,\alpha) = 2.
$$
Sixth, the resources are updated down $E=(0,0,6)$ and claimant 4 cannot be attended because the resources she claims are exhausted, as a result she gets 0. Therefore, we move to the next in the priority order. Again, claimant 6 cannot be attended because one of the resources she claims is exhausted, as a result she gets 0. So, we move to the next. Claimant 8 is attended:
$$
CSP_{8}^{\sigma}(\mathcal{I},N,E,c,\alpha) = 5.
$$
Seventh, the resources are updated down $E=(0,0,1)$ and the sequential procedure ends. The final allocation is
$$
CSP^{\sigma}(\mathcal{I},N,E,c,\alpha) = (3, 2, 4, 0, 5, 0, 3, 5).
$$
\end{example}

Once constrained sequential priority rules associated with an order have been defined, the constrained random arrival rule is simply defined as their average.

\begin{definition}\label{CRA}
Let $(\mathcal{I},N,E,c,\alpha) \in \mathcal{MBC}$, the \emph{constrained random arrival rule} for $(\mathcal{I},N,E,c,\alpha)$, $CRA(\mathcal{I},N,E,c,\alpha)$, is defined as follows:
$$
CRA_j(\mathcal{I},N,E,c,\alpha)=\frac{1}{n!}\sum_{\sigma\in\Sigma(N)}CSP_{j}^{\sigma}(\mathcal{I},N,E,c,\alpha),\:\forall j\in N,
$$
where $\Sigma(N)$ is the set of all possible orders of $N$.
\end{definition}

Although the definition of the rule is simple, it is immediately apparent that the biggest problem is in its computation. For Example \ref{ex1}, it would be necessary to calculate 8!=40320 sequential processes, which consumes a large amount of time.


\section{Properties}\label{properties}

In this section, we present several properties which are interesting in the context of MBC problems. These properties are related to efficiency, fairness, priority, monotonicity and consistency.

First, we introduce two concepts related to two claimants comparisons. In MBC situations claimants are characterized by two elements: their claims and the issues to which they claim. Therefore, both should be taken into account when establishing comparisons among them.

\begin{definition}
Let $(\mathcal{I},N,E,c,\alpha) \in \mathcal{MBC}$, and two claimants $j,k \in N$, we say they are {\em homologous}, if $\alpha(j) = \alpha(k)$; and they are called {\em equal}, if they are homologous and $c_j=c_k$.
\end{definition}

Next, we give a set of properties which are very natural and reasonable for an allocation rule in MBC situations.

The first property relates to efficiency. In allocation problems is desirable that the resources to be fully distributed, but in MBC situations this is not always possible (see Acosta-Vega, 2021). Therefore, a weaker version of that is considered in which only is required that there is no a feasible allocation in which at least one of the claimants receive more. This is established in the following axiom.

\begin{axiom}[PEFF]\label{PEFF}
Given a rule $R$, it satisfies {\em Pareto efficiency}, if for every problem $(\mathcal{I},N,E,c,\alpha) \in \mathcal{MBC}$, there is no a feasible allocation $a \in \mathbb{R}^N_{+}$ such that $a_j \geq R_j(\mathcal{I},N,E,c,\alpha), \forall j \in N$, with at least one strict inequality.
\end{axiom}

Note that $PEFF$ implies that at least the available amount of one issue is fully distributed, and no amount is left of an issue undistributed, if it is possible to do so. However, it does not require that all available amounts of the issues have to be fully distributed. On the other hand, a feasible allocation that satisfies the condition in Axiom \ref{PEFF} is called Pareto efficient. 

The second property states that equal claimants should receive the same in the final allocation. This is a basic requirement of fairness and non-arbitrariness. This is defined in the following axiom.

\begin{axiom}[ETE]\label{ETE}
Given a rule $R$, it satisfies {\em equal treatment of equals}, if for every problem $(\mathcal{I},N,E,c,\alpha)\in \mathcal{MBC}$ and every pair of equal claimants $j,k \in N$, $R_j(\mathcal{I},N,E,c,\alpha)=R_k(\mathcal{I},N,E,c,\alpha)$.
\end{axiom}

The third property is a requirement of robustness when some agents leave the problem with their allocations (see Thomson, 2011, 2018). In particular, when a subset of claimants leave the problem respecting the allocations to those who remain, then it seems reasonable that claimants who leave will receive the same in the new problem as they did in the original. This is formally formulated in the following axiom.

\begin{axiom}[CONS]\label{CONS}
Given a rule $R$, it satisfies {\em consistency}, if for every problem $(\mathcal{I},N,E,c,\alpha)\in \mathcal{MBC}$, and $N' \subset N$, it holds that 
$$
R_j(\mathcal{I},N,E,c,\alpha) = R_j(\mathcal{I}',N',E'^R,c|_{N'},\alpha), \text{ for all } j \in N',
$$
where $(\mathcal{I}',N',E'^R,c|_{N'},\alpha)\in \mathcal{MBC}$, called   the \emph{reduced problem associated with $N'$}, $\mathcal{I}'=\{i \in \mathcal{I}: \text{ there exists } k \in N' \text{such that }i \in \alpha(k)\}$, $E'^R=(e'^R_1,\ldots,e'^R_m)$ so that $e'^R_i=e_i- \sum_{j \in N\backslash N':i \in \alpha(j)}R_j(\mathcal{I},N,E,c,\alpha)$, for all $i \in \mathcal{I}'$, and $c|_{N'}$ is the vector whose coordinates correspond to the claimants in $N'$.
\end{axiom}

The fourth property is related to the priority order.

\begin{axiom}[PRI]\label{PRI}
Given a rule $R$, it satisfies {\em priority}, if for every problem $(\mathcal{I},N,E,c,\alpha)\in \mathcal{MBC}$ and every pair of homologous claimants $j,k \in N$, if $\sigma(j) < \sigma(k)$, $c_j -R_j(\mathcal{I},N,E,c,\alpha)\leq c_k -R_k(\mathcal{I},N,E,c,\alpha)$.
\end{axiom}

\begin{theorem}\label{CSPTh}
$CSP^{\sigma}$ for multi-issue bankruptcy problems with crossed claims satisfies $PEFF$, $CONS$, and $PRI$.
\begin{proof}We go axiom by axiom.
\begin{itemize}
\item $CSP^{\sigma}$ satisfies $PEFF$ and $PRI$ by definition.
\item Given $(\mathcal{I},N,E,c,\alpha)\in \mathcal{MBC}$ and $(\mathcal{I}',N',E'^{CPA},c|_{N'},\alpha)\in \mathcal{MBC}$ the reduced problem associated with $N' \subset N$. We assume that when some claimants leave the problem, they respect the priority order defined by $\sigma$ in the new problem. Therefore, if we follow the sequential procedure when a claimant in $N'$ has to be attended, what she can get is exactly the same as what she got in the original problem, since what the claimants who remained in the original problem keep is exactly what corresponded to them according to the order of priority defined by $\sigma$. Therefore, $CSP^{\sigma}$ satisfies $CONS$.
\end{itemize}
\end{proof}
\end{theorem}

It is well-known that the random arrival rule does not satisfy $CONS$, so neither does the constrained random arrival rule.

\begin{theorem}\label{CRATh}
$CRA$ for multi-issue bankruptcy problems with crossed claims satisfies $ETE$.
\begin{proof}The proof is straightforward by definition of $CRA$.
\end{proof}
\end{theorem}

In the context of multi-issue bankruptcy problems with crossed claims, the average operator does not preserve the property of Pareto efficiency and, therefore, CRA does not satisfy this basic property, as shown in the following proposition.

\begin{proposition}\label{PropnoPEFF}
$CRA$ rule does not satisfies $PEFF$.
\begin{proof}
Consider the problem with $M=\{1,2\}$, $N=\{1,2,3\}$, $E=(4,8)$, $c=(2,5,7)$, and $\alpha(1)=\{1\},\alpha(2)=\{1,2\},\alpha(3)=\{2\}$. $CSP$ for each order and $CRA$ are given by

\begin{table}[H]
\begin{center}
\begin{tabular}{cccc}
\textbf{Order / Claimant} & 1 & 2 & 3  \\ 
\toprule
123 & 2 & 2 & 6 \\
132 & 2 & 1 & 7 \\
213 & 0 & 4 & 4 \\
231 & 0 & 4 & 4 \\
312 & 2 & 1 & 7 \\
321 & 2 & 1 & 7 \\
\midrule
CRA & $\frac{8}{6}$ & $\frac{13}{6}$ & $\frac{35}{6}$ \\
\bottomrule
\end{tabular}
\end{center}
\caption{CSP and CRA for Proposition \ref{PropnoPEFF}.}\label{noPEFF}
\end{table}

Note that $\frac{8}{6} + \frac{13}{6} < 4$ and $\frac{8}{6} < 2$, therefore we can improve up to $\frac{11}{6}$ the allocation of agent 1 without exceeding any estate. Consequently, $CRA$ does not satisfy $PEFF$.
\end{proof}
\end{proposition}

A property that is satisfied by most of the rules for bankruptcy problems is resource monotonicity which is also used in the characterization of sequential priority rules (see Thomson, 2019, Th.11.11). This property simply says that if the available resource increases, allocations to claimants do not decrease. In the particular context of multi-issue bankruptcy problems with crossed claims, this property reads as follows.

\begin{axiom}[R-MON]\label{R-MON}
Given a rule $R$, it satisfies {\em resource monotonicity}, if for every problem $(\mathcal{I},N,E,c,\alpha)\in \mathcal{MBC}$ and $E' \geq E$, $R_j(\mathcal{I},N,E',c,\alpha) \geq R_j(\mathcal{I},N,E,c,\alpha)$ for all $j \in N$.
\end{axiom}

Although this property seems very weak, constrained sequential priority rules do not satisfy it as the following proposition shows.

\begin{proposition}
$CSP^{\sigma}$ for multi-issue bankruptcy problems with crossed claims does not satisfy $R-MON$.
\begin{proof}
Consider the problem with $\mathcal{I}=\{1,2,3\}$, $N=\{1,2,3,4,5,6,7,8\}$, $E=(9,12,7)$, $c=(3,5,4,3,5,4,4,5)$, and $\alpha(1)=\{1\},\alpha(2)=\{1\},\alpha(3)=\{1,2\},\alpha(4)=\{1,2\},\alpha(5)=\{2\},\alpha(6)=\{2,3\},\alpha(7)=\{2,3\},\alpha(8)=\{3\}$. If we take for example the priority order $\sigma = 1\:3\:5\:7\:2\:4\:6\:8$,
$$
CSP^{\sigma}(\mathcal{I},N,E,c,\alpha)=(3,2,4,0,5,0,3,4),
$$
and there is nothing left.

However, if $E'=(9,13,7)$, we have that
$$
CSP^{\sigma}(\mathcal{I},N,E',c,\alpha)=(3,2,4,0,5,0,4,3).
$$
Note that claimant 8 obtains less when the available resources increase. Therefore, $CSP^{\sigma}$ does not satisfy $R-MON$.
\end{proof}
\end{proposition}

Another monotonicity property is population monotonicity which says that if all claimants agree that a claimant $j$ will obtain her claim, then the remaining claimants should be worse off after claimant $j$ is fully compensated. This property is satisfied by many bankruptcy rules, including the random arrival rule. Moreover, the random arrival rule is characterized by using this property in Hwang and Wang (2009, Th.1). On the other hand, a property that appears recurrently related to the Shapley value and, therefore, to the random arrival rule is that of balanced contributions (Myerson, 1980; Hart and Mas-Colell, 1989). In the context of bankruptcy problems this property was used by Berganti\~nos and M\'endez-Naya (1997) to characterize the random arrival rule. Moreover, in the context of multi-issue bankruptcy problems was introduced by Lorenzo-Freire et al. (2007) and also used to characterize the random arrival rule. This property requires that claimant $j$ impacts to claimant $k$'s allocation what claimant $k$ impacts to claimant $j$'s allocation. Natural extensions of population monotonicity and balanced contributions to the context of multi-issue bankruptcy problems with crossed claims are the following.

\begin{axiom}[P-MON]{P-MON}
Given a rule $R$, it satisfies {\em population monotonicity}, if for every problem $(\mathcal{I},N,E,c,\alpha)$, and each $j \in N$,
$$
R_k(\mathcal{I},N,E,c,\alpha) \geq R_k(\mathcal{I}^{-j},N,E^{-j},c_{-j},\alpha), \text{ for all } k \in N\backslash {j},
$$
where $\mathcal{I}^{-j}$ is the set of issues for which claimants in $N\backslash \{j\}$ have claims; $E^{-j}=(e_1^{-j},\ldots,e_m^{-j})$ so that $e_i^{-j}=e_i - c_j$ if $i \in \alpha(j)$ and $e_i^{-j}=e_i$ otherwise; and $c_{-j}$ is the vector of claims from which $j-$th coordinate has been deleted.

\end{axiom}

\begin{axiom}[BAL]\label{BAL}
Given a rule $R$, it satisfies {\em balanced impact}, if for every problem $(\mathcal{I},N,E,c,\alpha)$, and every pair of claimants $j,k \in N$,
$$
R_j(\mathcal{I},N,E,c,\alpha) - R_j(\mathcal{I}^{-k},N,E^{-k},c_{-k},\alpha) = R_k(\mathcal{I},N,E,c,\alpha) - R_k(\mathcal{I}^{-j},N,E^{-j},c_{-j},\alpha),
$$
where $\mathcal{I}^{-h}$ is the set of issues for which claimants in $N\backslash \{h\}$ have claims; $E^{-h}=(e_1^{-h},\ldots,e_m^{-h})$ so that $e_i^{-h}=e_i - c_h$ if $i \in \alpha(h)$ and $e_i^{-h}=e_i$ otherwise; and $c_{-h}$ is the vector of claims from which $h-$th coordinate has been deleted.
\end{axiom}

As mentioned above, this type of properties are used to characterize Shapley-like solutions as the random arrival is, but in the context of multi-issue bankruptcy problems with crossed claims, the constrained random arrival rule does not satisfy them as the following proposition shows.

\begin{proposition}\label{propbal}
$CRA$ rule satisfies neither $P-MON$ nor $BAL$.
\begin{proof}
Consider the situation described by Figure \ref{figex2} with three issues $\mathcal{I}=\{1,2,3\}$, and three claimants $N=\{1,2,3\}$, the available resources are given by $E=(4,5,7)$, the claims are given by the vector $c=(3,4,5)$, and finally, $\alpha(1)=\{1,2\},\alpha(2)=\{2,3\},\alpha(3)=\{3\}$.

\begin{figure}
\begin{center}
\begin{tikzpicture}
\node[draw] (c3) at (8,0) {{\small $c_3=5$}};
\node[draw] (c2) at (6,0) {{\small $c_2=4$}};
\node[draw] (c1) at (1.5,0) {{\small $c_1=3$}};

\node[draw] (e1) at (1.5,2.5) {{\small $e_1=4$}};
\node[draw] (e2) at (5,2.5) {{\small $e_2=5$}};
\node[draw] (e3) at (8,2.5) {{\small $e_3=7$}};

\draw [->] (c1) -- (e1);

\draw [->] (c1) -- (e2);
\draw [->] (c2) -- (e2);

\draw [->] (c2) -- (e3);
\draw [->] (c3) -- (e3);

\end{tikzpicture}
\end{center}
\caption{Multi-issue problem with crossed claims with three issues and three claimants.}\label{figex2}
\end{figure}
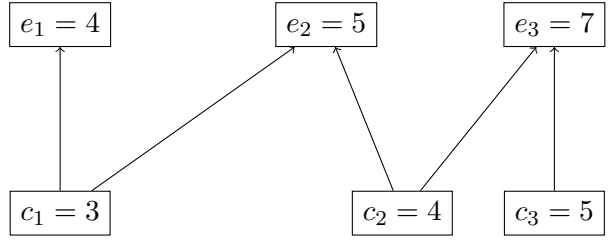

In order to calculate $CRA$, we need to consider all possible orders of arrival. These calculations are summarized in Table \ref{todos} below. 

\begin{table}[H]
\begin{center}
\begin{tabular}{cccc}
\textbf{Order / Claimant} & 1 & 2 & 3  \\ 
\toprule
123 & 3 & 2 & 5  \\
132 & 3 & 2 & 5 \\ 
213 & 1 & 4 & 3 \\
231 & 1 & 4 & 3 \\
312 & 3 & 2 & 5 \\
321 & 3 & 2 & 5 \\
\midrule
CRA & $\frac{7}{3}$ & $\frac{8}{3}$ & $\frac{13}{3}$\\ 
\bottomrule
\end{tabular}
\end{center}
\caption{Constrained random arrival rule for Proposition \ref{propbal}.}\label{todos}
\end{table}

In Table \ref{parados} the calculations to obtain $CRA$ for the two claimants subproblems are shown.

\begin{table}[H]
\begin{center}
\begin{tabular}{ccccccccc}
\textbf{Order} & 1 & 2 & \textbf{Order} & 1 & 3 & \textbf{Order} & 2 & 3  \\ 
\toprule
12 & 3 & 2 & 13 & 1 & 3 & 23 & 2 & 5  \\
21 & 3 & 2 & 31 & 1 & 3 & 32 & 2 & 5 \\ 
\midrule
CRA & 3 & 2 &   & 1 & 3 & & 2 & 5 \\ 
\bottomrule
\end{tabular}
\end{center}
\caption{Constrained random arrival rule for the two claimants subproblems in Proposition \ref{propbal}.}\label{parados}
\end{table}

Now, we observe that $CRA_1(\mathcal{I},N,E,c,\alpha) =\frac{7}{3} < 3 = CRA_1(\mathcal{I}^{-3},N^{-3},E^{-3},c_{-3},\alpha)$. Therefore, claimant 1 is better off after claimant 3 leaves with her claim. Consequently, $CRA$ does not satisfy $P-MON$.

On the other hand, we have that
$$
CRA_1(\mathcal{I},N,E,c,\alpha) - CRA_1(\mathcal{I}^{-2},N^{-2},E^{-2},c_{-2},\alpha) =\frac{7}{3} -1 = \frac{4}{3},
$$
and
$$
CRA_2(\mathcal{I},N,E,c,\alpha) - CRA_2(\mathcal{I}^{-1},N^{-1},E^{-1},c_{-1},\alpha) =\frac{8}{3} -2 = \frac{2}{3}.
$$
Therefore, $CRA$ does not satisfy $BAL$.
\end{proof}
\end{proposition}

The previous results show the complexity of finding properties that allow axiomatic characterizations of sequential priority rules and the random arrival rule in the context of multi-issue bankruptcy problems with crossed claims, so it is necessary to look for perhaps more specific properties (and likely more technical) to achieve it.


\section{An alternative extension of the random arrival rule for MBC}

In this section, we introduce a modification of the constrained random arrival rule for multi-issue bankruptcy problems with crossed claims, that we call CRA*. In particular, the CRA* consists of two levels of arrival orders. However, for this rule to be well defined, it is first necessary that the claims are truncated by the minimum of their associated endowments to avoid inconsistencies, i.e., we must take as claims $c'_i=\min\{c_i,\min_{k \in \alpha(i)}\{e_k\}\}, i \in N$, and this must be done before calculating each RA with the updated estates and claims. Once the claims have been truncated, first, all the possible orders of the issues are taken. For the first issue in the order, we calculate the RA rule, then the estates and claims are updated down and the RA rule is calculated for the second issue in the order, and so on until the last issue in the order. Once we have obtained an allocation for each possible order of issues we take their average. Note that CRA* can be seen as a new and different allocation from CRA, but we use the same rationale behind. Moreover, it is also an extension of the random arrival rule since for one issue bankruptcy problems both coincide. The following example illustrates how CRA* works.

\begin{example}\label{exCRA*}
Consider the situation described by Figure \ref{FigexCRA*} with $\mathcal{I}=\{1,2,3\}$, $N=\{1,2,3,4,5\}$, $E=(9,10,8)$, $c=(3,4,3,6,5)$, and $\alpha(1)=\{1\},\alpha(2)=\{1,2\},\alpha(3)=\{1,2\},\alpha(4)=\{2,3\},\alpha(5)=\{3\}$.

\begin{figure}
\begin{center}
\begin{tikzpicture}
\node[draw] (c5) at (8,0) {{\small C5}};
\node[draw] (c4) at (6,0) {{\small C4}};

\node[draw] (c3) at (4.5,0) {{\small C3}};
\node[draw] (c2) at (3,0) {{\small C2}};
\node[draw] (c1) at (1.5,0) {{\small C1}};

\node[draw] (e1) at (1.5,2.5) {{\small E1}};
\node[draw] (e2) at (5,2.5) {{\small E2}};
\node[draw] (e3) at (8,2.5) {{\small E3}};

\draw [->] (c1) -- (e1);
\draw [->] (c2) -- (e1);
\draw [->] (c3) -- (e1);

\draw [->] (c2) -- (e2);
\draw [->] (c3) -- (e2);
\draw [->] (c4) -- (e2);

\draw [->] (c4) -- (e3);
\draw [->] (c5) -- (e3);

\end{tikzpicture}
\end{center}
\caption{Multi-issue bankruptcy problem with crossed claims in Example \ref{exCRA*}.}\label{FigexCRA*}
\end{figure}
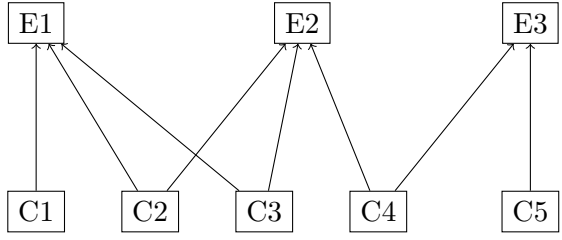

\begin{itemize}
\item Consider the order of issues {\bf 123}.

We start with issue 1. The bankruptcy problem associated with this issue is $e_1=9, c_1=3, c_2=4, c_3=3$. No need to truncate claims.

\begin{table}[H]
\begin{center}
\begin{tabular}{cccc}
\textbf{Order / Claimant} & 1 & 2 & 3  \\ 
\toprule
123 & 3 & 4 & 2  \\
132 & 3 & 3 & 3 \\ 
213 & 3 & 4 & 2 \\
231 & 2 & 4 & 3 \\
312 & 3 & 3 & 3 \\
321 & 2 & 4 & 3 \\
\midrule
RA & $2\frac{2}{3}$ & $3\frac{2}{3}$ & $2\frac{2}{3}$\\ 
\bottomrule
\end{tabular}
\end{center}
\end{table}

UPDATE PHASE:
\begin{enumerate}
\item We update the estates and claims: $e'_1=0, e'_2=3\frac{2}{3}, e'_3=8, c'_1=\frac{1}{3}, c'_2=\frac{1}{3}, c'_3=\frac{1}{3}, c'_4=6, c'_5=5$.
\item We truncate the claims by the estates: $c''_1=0, c''_2=0, c''_3=0, c''_4=3\frac{2}{3}, c''_5=5$.
\end{enumerate}

Now we continue with issue 2. The bankruptcy problem associated with this issue is $e'_2=3\frac{2}{3}, c''_2=0, c''_3=0, c''_4=3\frac{2}{3}$. Therefore, $RA = (0,0,3\frac{2}{3})$.

\bigskip

UPDATE PHASE:
\begin{enumerate}
\item We update the estates and claims: $e'_1=0, e'_2=0, e'_3=4\frac{1}{3}, c'_1=0, c'_2=0, c'_3=0, c'_4=0, c'_5=5$.
\item We truncate the claims by the estates: $c''_1=0, c''_2=0, c''_3=0, c''_4=0, c''_5=4\frac{1}{3}$.
\end{enumerate}

We end with issue 3. The bankruptcy problem associated with this issue is $e'_3=4\frac{1}{3}, c''_4=0, c''_5=4\frac{1}{3}$. Therefore, $RA = (0,4\frac{1}{3})$.

Therefore, the allocation for the order of issues 123 is $(2\frac{2}{3}, 3\frac{2}{3}, 2\frac{2}{3}, 3\frac{2}{3}, 4\frac{1}{3})$.

\bigskip

\item Consider the order of issues {\bf 132}.

We start with issue 1. This case is identical to order 123. Therefore, $RA=(2\frac{2}{3}, 3\frac{2}{3}, 2\frac{2}{3})$. The update phase is also the same, but we write it down for the sake of completeness.

\bigskip

UPDATE PHASE:
\begin{enumerate}
\item We update the estates and claims: $e'_1=0, e'_2=3\frac{2}{3}, e'_3=8, c'_1=\frac{1}{3}, c'_2=\frac{1}{3}, c'_3=\frac{1}{3}, c'_4=6, c'_5=5$.
\item We truncate the claims by the estates: $c''_1=0, c''_2=0, c''_3=0, c''_4=3\frac{2}{3}, c''_5=5$.
\end{enumerate}

Now we continue with issue 3. The bankruptcy problem associated with this issue is $e'_3=8, c''_4=3\frac{2}{3}, c''_5=5$. Therefore, $RA = (3\frac{1}{3},4\frac{2}{3})$.

\bigskip

UPDATE PHASE:
\begin{enumerate}
\item We update the estates and claims: $e'_1=0, e'_2=\frac{1}{3}, e'_3=0, c'_1=0, c'_2=0, c'_3=0, c'_4=\frac{1}{3}, c'_5=\frac{1}{3}$.
\item We truncate the claims by the estates: $c''_1=0, c''_2=0, c''_3=0, c''_4=0, c''_5=0$.
\end{enumerate}

We end with issue 2, but there is no bankruptcy problem since $c''_2=0, c''_3=0, c''_4=0$. 

Therefore, the allocation for the order of issues 132 is $(2\frac{2}{3}, 3\frac{2}{3}, 2\frac{2}{3}, 3\frac{1}{3}, 4\frac{2}{3})$.

\item Consider the order of issues {\bf 213}.

We start with issue 2. The bankruptcy problem associated with this issue is $e_2=10, c_2=4, c_3=3, c_4=6$. No need to truncate claims.

\begin{table}[H]
\begin{center}
\begin{tabular}{cccc}
\textbf{Order / Claimant} & 2 & 3 & 4  \\ 
\toprule
234 & 4 & 3 & 3  \\
243 & 4 & 0 & 6 \\ 
324 & 4 & 3 & 3 \\
342 & 1 & 3 & 6 \\
423 & 4 & 0 & 6 \\
432 & 1 & 3 & 6 \\
\midrule
RA & $3$ & $2$ & $5$\\ 
\bottomrule
\end{tabular}
\end{center}
\end{table}

UPDATE PHASE:
\begin{enumerate}
\item We update the estates and claims: $e'_1=4, e'_2=0, e'_3=3, c'_1=3, c'_2=1, c'_3=1, c'_4=1, c'_5=5$.
\item We truncate the claims by the estates: $c''_1=3, c''_2=0, c''_3=0, c''_4=0, c''_5=3$.
\end{enumerate}

Now we continue with issue 1 and then issue 3, but in these cases it is immediate that $RA=(3,0,0)$ and $RA=(0,3)$.

Therefore, the allocation for the order of issues 213 is $(3,3,2,5,3)$.

\item Consider the order of issues {\bf 231}. This case is completely identical to the order of issues 213, therefore, the allocation for the order of issues 231 is $(3,3,2,5,3)$.

\item Consider the order of issues {\bf 312}.

We start with issue 3. The bankruptcy problem associated with this issue is $e_3=10, c_4=6, c_5=5$. No need to truncate claims.

\begin{table}[H]
\begin{center}
\begin{tabular}{ccc}
\textbf{Order / Claimant} & 4 & 5  \\ 
\toprule
45 & 6 & 2  \\
54 & 3 & 5  \\ 
\midrule
RA & $4.5$ & $3.5$ \\ 
\bottomrule
\end{tabular}
\end{center}
\end{table}

UPDATE PHASE:
\begin{enumerate}
\item We update the estates and claims: $e'_1=9, e'_2=5.5, e'_3=0, c'_1=3, c'_2=4, c'_3=3, c'_4=1.5, c'_5=1.5$.
\item We truncate the claims by the estates: $c''_1=3, c''_2=4, c''_3=3, c''_4=0, c''_5=0$.
\end{enumerate}

Now we continue with issue 1. The bankruptcy problem associated with this issue is $e'_1=9, c''_1=3, c''_2=4, c''_3=3$.

\begin{table}[H]
\begin{center}
\begin{tabular}{cccc}
\textbf{Order / Claimant} & 1 & 2 & 3  \\ 
\toprule
123 & 3 & 4 & $1.5^a$  \\
132 & 3 & $2.5^b$ & 3 \\ 
213 & 3 & 4 & $1.5^a$ \\
231 & $3^c$ & 4 & $1.5^a$ \\
312 & 3 & $2.5^b$ & 3 \\
321 & $3^c$ & $2.5^b$ & 3 \\
\midrule
RA & $3$ & $3.25$ & $2.25$\\ 
\bottomrule
\end{tabular}
\end{center}
\caption{(a) Claimant 3 cannot receive 2 units because the sum of the allocations of claimants 2 and 3 would exceed $e'_2=5.5$. (b) Claimant 2 cannot receive 3 units because the sum of the allocations of claimants 2 and 3 would exceed $e'_2=5.5$. (c) Claimant 1 can receive more than 2 because by receiving claimants 2 and 3 less she can receive more.}
\end{table}

UPDATE PHASE:
\begin{enumerate}
\item We update the estates and claims: $e'_1=0.5, e'_2=0, e'_3=0, c'_1=0, c'_2=0.75, c'_3=0.75, c'_4=0, c'_5=0$.
\item We truncate the claims by the estates: $c''_1=0, c''_2=0, c''_3=0, c''_4=0, c''_5=0$.
\end{enumerate}

We end with issue 2, but there is no bankruptcy problem.

Therefore, the allocation for the order of issues 312 is $(3,3.25,2.25,4.5,3.5)$.

\item Consider the order of issues {\bf 321}.

We start with issue 3. The situation is the same as the order 312. Therefore, $RA=(4.5,3.5)$.

\bigskip

UPDATE PHASE:
\begin{enumerate}
\item We update the estates and claims: $e'_1=9, e'_2=5.5, e'_3=0, c'_1=3, c'_2=4, c'_3=3, c'_4=1.5, c'_5=1.5$.
\item We truncate the claims by the estates: $c''_1=3, c''_2=4, c''_3=3, c''_4=0, c''_5=0$.
\end{enumerate}

Now we continue with issue 2. The bankruptcy problem associated with this issue is $e'_2=5.5, c''_2=4, c''_3=3, c''_4=0$. $RA=(3.25,2.25,0)$.

\bigskip

UPDATE PHASE:
\begin{enumerate}
\item We update the estates and claims: $e'_1=3.5, e'_2=0, e'_3=0, c'_1=3, c'_2=0.75, c'_3=0.75 c'_4=0, c'_5=0$.
\item We truncate the claims by the estates: $c''_1=3, c''_2=0, c''_3=0, c''_4=0, c''_5=0$.
\end{enumerate}

We end with issue 1, but there is no bankruptcy problem and $RA=(3,0,0)$.

Therefore, the allocation for the order of issues 321 is $(3,3.25,2.25,4.5,3.5)$.

\end{itemize}

Finally, CRA* is determined by calculating the average of all the allocations obtained.

\begin{table}[H]
\begin{center}
\begin{tabular}{cccccc}
\textbf{Order / Claimant} & 1 & 2 & 3 & 4 & 5 \\ 
\toprule
123 & $2\frac{2}{3}$ & $3\frac{2}{3}$ & $2\frac{2}{3}$ & $3\frac{2}{3}$ & $4\frac{1}{3}$  \\
132 & $2\frac{2}{3}$ & $3\frac{2}{3}$ & $2\frac{2}{3}$ & $3\frac{1}{3}$ & $4\frac{2}{3}$ \\ 
213 & 3 & 3 & 2 & 5 & 3 \\
231 & 3 & 3 & 2 & 5 & 3 \\
312 & 3 & 3.25 & 2.25 & 4.5 & 3.5 \\
321 & 3 & 3.25 & 2.25 & 4.5 & 3.5 \\
\midrule
CRA* & $2\frac{8}{9}$ & $3\frac{11}{36}$ & $2\frac{11}{36}$ & $4\frac{1}{3}$ & $3\frac{2}{3}$\\ 
\bottomrule
\end{tabular}
\end{center}
\end{table}

Note that issue 3 is completely distributed, but issues 1 and 2 are not. Moreover, it is easy to check that the allocation given by CRA* is not Pareto efficient. 

\end{example}

Formally, CRA* is defined as follows.

\begin{definition}\label{CSP}
Let $(\mathcal{I},N,E,c,\alpha) \in \mathcal{MBC}$, $CRA*(\mathcal{I},N,E,c,\alpha)$ is defined as follows:

Let $\omega : \{1,2,\ldots,|\mathcal{I}|\} \rightarrow \mathcal{I}$ a one-to-one application, and $\Omega(\mathcal{I})$ the set of all those applications. Let $N^i=\{j \in N: i \in \alpha(j)\}$. Let $\xi^i : N^i \rightarrow \{1,2,\ldots,|N^i|\}$ a one-to-one application, and $\Xi(N^i)$ the set of all those applications.

\begin{itemize}

\item For each $\omega \in \Omega(\mathcal{I})$:

\begin{itemize}
\item Initialization:
\begin{enumerate} 
\item $e_i^{U(0)}=e_i, i \in \mathcal{I}$, and $c_j^{U(0)}=c_j, j \in N$. 
\item $c_j^{UT(0)}=\min\left\{c_j^{U(0)},\min\{e_i^{U(0)}: i \in \alpha(j)\}\right\}, j \in N$.
\end{enumerate}
\item From $k=1$ to $|\mathcal{I}|$, do
\begin{enumerate}
\item For $\omega(k)$, we consider the bankrupty problem $(N^{\omega(k)}, e_{\omega(k)}^{U(k-1)},\{c_j^{UT(k-1)}\}_{j \in N^{\omega(k)}})$.
\item We calculate $m_j^{\xi^{\omega(k)}}, j \in N^{\omega(k)}$ as follows:
$$
m_j^{\xi^{\omega(k)}}=min\left\{ c_j^{UT(k-1)},\:max\left\{0,\:\min_{i \in \alpha(j)}\{e_i^{U(k-1)}-\sum_{\substack{
h \in N^{\omega(k)}:i \in \alpha(h) \\[0.1cm]
\xi^{\omega(k)}(h)<\xi^{\omega(k)}(j)}
}
c_h^{UT(k-1)}\}\right\}\right\}.
$$
\item We calculate
$$RA^{\omega(k)}_j=\frac{1}{|N^{\omega(k)}|!}\sum_{\xi^{\omega(k)} \in \Xi(N^{\omega(k)})}m_j^{\xi^{\omega(k)}}, j \in N^{\omega(k)}.$$

\item UPDATE PHASE
\begin{enumerate}
\item $e_i^{U(k)}=e_i^{U(k-1)} - \sum_{j \in N^{\omega(k)}:i \alpha(j)} RA^{\omega(k)}_j, i \in \mathcal{I}$; and $c_j^{U(k)}=c_j^{UT(k-1)} - RA^{\omega(k)}_j, j \in N^{\omega(k)}; c_j^{U(k)}=c_j^{UT(k-1)}, j \in N \backslash N^{\omega(k)}$. 
\item $c_j^{UT(k)}=\min\{c_j^{U(k)},\min\{e_i^{U(k)}: i \in \alpha(j)\}, j \in N$.
\end{enumerate}
\end{enumerate}

\item The last step is to calculate the allocation associated with $\omega$:
$$
A_j^{\omega}=\sum_{k =1}^{|\mathcal{I}|}\delta(j,\omega(k))RA^{\omega(k)}_j, j \in N,
$$
where $\delta(j,\omega(k)) = 1$ if $\omega(k) \in \alpha(j)$ and $0$ otherwise.
\end{itemize}

\item Finally, we calculate CRA* for each $j \in N$ as follows:
$$
CRA_j^{*}=\frac{1}{|\mathcal{I}|!}\sum_{\omega \in \Omega(\mathcal(I))}A_j^{\omega}.
$$
\end{itemize}
\end{definition}

The following examples shows that CRA and CRA* do not coincide in general.

\begin{example}
Consider the multi-issue bankruptcy problem with crossed claims in Proposition \ref{propbal}. The allocation for each order of issues and CRA* are shown in the following table:

\begin{table}[H]
\begin{center}
\begin{tabular}{cccc}
\textbf{Order of issues / Claimant} & 1 & 2 & 3  \\ 
\toprule
123 & 3 & 2 & 5  \\
132 & 3 & 2 & 5 \\ 
213 & 2 & 3 & 4 \\
231 & 2 & 3 & 4 \\
312 & 2 & 3 & 4 \\
321 & 2 & 3 & 4 \\
\midrule
CRA* & $\frac{7}{3}$ & $\frac{8}{3}$ & $\frac{13}{3}$\\ 
\bottomrule
\end{tabular}
\end{center}
\end{table}

We observe that CRA and CRA* coincide. Moreover, if claimant 3 leaves the problem with her claim, CRA* of the new problem is $(3,2)$, therefore, CRA* does not satisfies P-MON. Furthermore, if claimant 1 leaves the problem with her claim, CRA* of the new problem is $(2,5)$; and if claimant 2 leaves the problem with her claim, CRA* of the new problem is $(1,3)$, therefore CRA* does not satisfies BAL.
\end{example}

\begin{example}
Consider the multi-issue bankruptcy problem with crossed claims in Proposition \ref{PropnoPEFF}. The allocation for each order of issues and CRA* are shown in the following table:

\begin{table}[H]
\begin{center}
\begin{tabular}{cccc}
\textbf{Order of issues / Claimant} & 1 & 2 & 3  \\ 
\toprule
12 & 1 & 3 & 5  \\
21 & 1.5 & 2.5 & 5.5 \\ 
\midrule
CRA* & $1.25$ & $2.75$ & $5.25$\\ 
\bottomrule
\end{tabular}
\end{center}
\end{table}

We now observe that CRA and CRA* do not coincide. Moreover, for this example CRA* is Pareto efficient, but CRA does not.
\end{example}

Finally, we observe that we have the same problems with the properties as for CRA, therefore, it is not easy to find a characterization of CRA* using the usual properties in the characterizations of the random arrival rule for bankruptcy problems.


\section{Conclusions}
\label{conc}

In many allocation problems the concept of priority is relevant, for example, in the legislation related to the liquidation of a company through bankruptcy in many countries, an order of priority is established to satisfy the claims. First, wage claims are settled, then taxes are paid. Creditor claims are then satisfied, and finally shareholder claims are addressed. Therefore, sequential priority rules, although simple, are not strange in real life. In this work an extension of the sequential priority rules has been introduced in the context of multi-issue bankruptcy problems with crossed claims. Next, by means of the average of all these rules, the extension of the random arrival rule to this new context is defined.

In the analysis of the properties satisfied by the sequential priority rules and the random arrival rule, it has been shown that the natural extensions of the properties used in the characterization of these rules in the context of bankruptcy problems and multi-issue bankruptcy problems are not satisfied by the constrained sequential priority rules and the constrained random arrival rule. Moreover, CRA* does not satisfies either the usual properties in the characterization of the random arrival rule. Therefore, the characterizations of the rules introduced in this work will require properties that are perhaps too technical or very \textit{ad hoc}, which could detract from a simple interpretation.

The reason that these properties are not satisfied is that there is no efficiency in the context of multi-issue problems with crossed claims but rather a weaker concept such as Pareto efficiency. This means that the total quantity distributed can change from one problem to another, so the results from one problem to another are not easily comparable. This is a problem but it also indicates us that these problems are of theoretical interest because they are not a mere and simple extension of bankruptcy problems but rather have a different structure that makes them interesting for further study.

\section*{Acknowledgment}
This work is part of the R\&D\&I project grant PGC2018-097965-B-I00, funded by MCIN/ AEI/10.13039/501100011033/ and by "ERDF A way of making Europe"/EU. The authors are grateful for this financial support. Joaqu\'in S\'anchez-Soriano also acknowledges financial support from the Generalitat Valenciana under the project PROMETEO/2021/063.

\end{document}